\documentclass[a4paper,10pt,reqno]{amsart}
\usepackage{verbatim}
\usepackage{amssymb}      
       
\usepackage{enumerate}       
\usepackage[active]{srcltx}       
\numberwithin{equation}{section}       
\usepackage{paralist} 
\usepackage{minibox}
       
\usepackage{t1enc}       
\usepackage[utf8x]{inputenc}       
\usepackage{mathabx}
 \usepackage{tikz}    
 \usetikzlibrary{shapes}

 \newcommand{\snode}[2]{\node (#1) [ pro] {\tiny #2};}
       
\newtheorem{theorem}{Theorem}[section]        
\newtheorem{lemma}[theorem]{Lemma}

\newtheorem{problem}[theorem]{Problem}

\newtheorem{corollary}[theorem]{Corollary}       
\newtheorem{claim}{Claim}[theorem]

\theoremstyle{definition}       
\newtheorem{definition}[theorem]{Definition}       
\theoremstyle{remark}       
         
\newcommand{\mc}[1]{\mathcal{#1}}       
\newcommand{\mbb}[1]{\mathbb{#1}}

\newcommand{\setm}{\setminus}       
\newcommand{\empt}{\emptyset}       
\newcommand{\subs}{\subset}       
\newcommand{\dom}{\operatorname{dom}}       
\newcommand{\ran}{\operatorname{ran}}

\def\<{\left\langle}       
\def\>{\right\rangle}       
       
\makeatletter      
  \@namedef{subjclassname@2020}{%
  \textup{2020} Mathematics Subject Classification}    
   \makeatother    
       
\author[T. Csern\'ak]{
Tam\'as Csern\'ak }
\address{Eötvös University of Budapest, Hungary}
\email{tamas@csernak.com}

\author[L. Soukup]{
Lajos Soukup }
\address{Alfr{\'e}d R{\'e}nyi Institute of Mathematics }
\email{soukup@renyi.hu}

\title[Infinite Combinatorics revisited]{Infinite Combinatorics revisited in the absence of Axiom of Choice}       
\thanks{The preparation of this paper was supported
by NKFIH grant no. K129211.}
\date{\today}       
       

\begin{document}

\newcommand{\dsystem}[3]{{[#1]}^{<#2}\to \Delta(#3)}
\newcommand{\freesetsmall}[3]{Free(#1,#2,#3)}
\newcommand{\freeset}[3]{Free(#1\mapsto {[#1]}^{<#2}, #3)}
\newcommand{\freeUnionsmall}[3]{FreeUnion(#1,#2,#3)}
\newcommand{\freeUnion}[3]{FreeUnion(#1\mapsto {[#1]}^{<#2}, #3)}

\newcommand{\KK}[2]{#1[#2]}
\newcommand{\KKG}[1]{\KK{#1}{\mc G}}
\newcommand{\KKGU}[1]{\KKG{\underline{#1}}}

\begin{abstract}
    We investigate the provability of classical combinatorial theorems in ZF.
    Using purely combinatorial arguments, we establish the following results for each
infinite cardinal ${\kappa}\in On$, 
    \begin{inparaenum}[(1)]
\item ${\kappa}^+\to ({\kappa},{\omega}+1)$, 
\item any family $\mc A\subs [{On}]^{<{\omega}}$ of size ${\kappa}^+$
contains a $\Delta$-system of size ${\kappa}$,
\item any regressive function $f:{\kappa}^+\to {\kappa}^+$ is constant on a set
of size ${\kappa}$,
\item given a set mapping $F:{\kappa}\to {[{\kappa}]}^{<{\omega}}$, 
the set  ${\kappa}$ has a  partition into ${\omega}$-many $F$-free sets,  
\item there is a cardinal ${\lambda}\in On$
such that ${\lambda}\to ({\kappa}^+)^2_{\kappa}$.
\end{inparaenum}

By employing Karagila's method of absoluteness, we demonstrate the following results for each uncountable cardinal 
${\kappa}\in On$,
\begin{inparaenum}[(1)]
        \addtocounter{enumi}{5}
\item given a set mapping $F:{\kappa}\to {[{\kappa}^]}^{<{\omega}}$,  there is 
an $F$-free set of cardinality ${\kappa}$, 
\item for each natural number $n$, every family $\mc A\subs {[{\kappa}]}^{{\omega}}$
with $|A\cap B|\le n$ for $\{A,B\}\in {[\mc A]}^{2}$ has property $B$,  
\end{inparaenum}

In contrast to statement (7), we show that the following statement is not provable from ZF +
 $cf({\omega}_1)={\omega}_1$:
\begin{enumerate}[($8^*$)]
\item every family $\mc A\subs {[{\omega}_1]}^{{\omega}}$
with $|A\cap B|\le 1$ for $\{A,B\}\in {[\mc A]}^{2}$ is {\em essentially disjoint}
(i.e. there is a 
        function $f$ with $\dom(f)=\mc A$ and $f(A)\in [A]^{<{\omega}}$
        such that $\{A\setm f(A):A\in \mc A\}$ is a family of pairwise disjoint sets). 
\end{enumerate}

The following statements are not provable in ZF, but they are 
equivalent within ZF:
    \begin{inparaenum}[(i)]
     \item $cf({\omega}_1)={\omega}_1$, 
     \item ${\omega}_1\to ({\omega}_1,{\omega}+1)^2$,
     \item  any family $\mc A\subs [{On}]^{<{\omega}}$ of size ${\omega}_1$
     contains a $\Delta$-system of size ${\omega}_1$.
    \end{inparaenum}

We say that a function $f$ is a {\em uniform denumeration on ${\omega}_1$} iff 
$\dom(f)={\omega}_1$ and for every ${\alpha}<{\omega}_1$, $f({\alpha})$ is a function from 
${\omega}$ onto ${\alpha}$.

It is evident that 
the existence of a uniform denumeration of ${\omega}_1$ implies 
$cf({\omega}_1)={\omega}_1$. We prove that the failure of the reverse implication 
is equiconsistent with the existence of an inaccessible cardinal.  

\end{abstract}
\maketitle

\section{Introduction}
The classical theorems of infinite combinatorics, such as the $\Delta$-system theorem, the pressing down lemma, the Erdős-Dushnik-Miller partition theorem, and the free set theorem of Hajnal, have traditionally relied heavily on the Axiom of Choice. Our investigation began with the question of whether these theorems can be proven without invoking the Axiom of Choice.

In \cite{Ka14}, Karagila demonstrated that the Erdős-Dushnik-Miller partition theorem ${\kappa}\to ({\kappa},{\omega})$ can be derived from ZF without using the Axiom of Choice. However, his argument is not combinatorial; it is based on the concept of absoluteness between different ZF models. In this paper, our aim is to utilize purely combinatorial arguments whenever possible, in order to present our results in a manner that is accessible and appealing to a wider audience.

\pagebreak
\subsection*{Definitions and outline of our  results.} \nopagebreak 
\begin{definition}
If ${\mu},{\lambda}\in On $ are cardinals,
write $Regressive({\mu},{\lambda})$ iff for 
each a regressive function $f:{\mu}\to {\mu}$  
there is ${\alpha}<{\mu}$ such that 
$|f^{-1}\{{\alpha}\}|\ge {\lambda}$.      
\end{definition}

In \cite{Du31-reg} Dushnik proved that $Regressive({\kappa}^+,{\kappa}^+)$ holds for each infinite cardinal ${\kappa}$
in ZFC.

\begin{definition}
If ${\mu},{\lambda}$ are infinite cardinals, write
$$\dsystem{{\mu}}{{\omega}}{{\lambda}}$$
iff any family $\mc A$ of finite sets of ordinals with $|\mc A|={\mu}$ contains
a $\Delta$-system $\mc B$ with $|\mc B|={\lambda}$.
\end{definition}

Shanin proved in \cite{Sa46} that 
 $$\dsystem{{\mu}}{{\omega}}{\mu}$$ holds for each uncountable  regular cardinal 
${\mu}$ in ZFC (see \cite[Theorem 9.18]{Je82}).

\subsection*{Equivalences and weaker statements in ZF}

ZF is not enough to prove that $cf({\kappa}^+)={\kappa}^+$, $Regressive({\kappa}^+,{\kappa}^+)$
or $\dsystem{{\kappa}^+}{{\omega}}{{\kappa}^+}$, where ${\kappa}\in On $ is an infinite cardinal,
but in Theorem \ref{tm:equivalence-in-ZF} we show that these statements are equivalent to each other 
in ZF, and we can get the following corollary:

\begin{corollary}\label{cr:eq}[ZF]
The following statements are equivalent:
\begin{inparaenum}[\itshape a\upshape)]
    \item $cf({\omega}_1)={\omega}_1$, 
    \item ${\omega}_1\to ({\omega}_1,{\omega}+1)$,
    \item $Regressive({\omega}_1,{\omega}_1)$,
    \item $\dsystem{{\omega}_1}{{\omega}}{{\omega}_1}$.
    \end{inparaenum}
\end{corollary}

Moreover, in Theorem \ref{tm:elementary-in-ZF} we also show that certain mild weakening 
of these ZFC results can be proved in ZF, in particular,  
\begin{corollary}\label{crweak-:true}[ZF]
    \begin{inparaenum}[\itshape a\upshape)]
        \item ${\omega}_2\to ({\omega}_1)^1_{{\omega}_1}$, 
        \item ${\omega}_2\to ({\omega}_1,{\omega}+1)$,
        \item $Regressive({\omega}_2,{\omega}_1)$,
        \item $\dsystem{{\omega}_2}{{\omega}}{{\omega}_1}$.
        \end{inparaenum}
    
    \end{corollary}

Next we present  some classical ZFC results which are provable in ZF.

\begin{definition}
    Given a set mapping  $F:{\kappa}\to {[{\kappa}]}^{<{\mu}}$ we say that a set 
    $X\subs {\kappa}$ is {\em $F$-free} iff ${\alpha}\notin F({\beta})$ for each $\{{\alpha},{\beta}\}\in {[X]}^{2}$. 
            If ${\kappa},{\lambda},{\mu}$ are cardinals, write 
           \begin{enumerate}[(i)]
       \item         $\freeset{{{\kappa}}}{{\mu}}{{\lambda}}$ iff for each
       set-mapping $F:{\kappa}\to {[{\kappa}]}^{<{\mu}}$ there is an $F$-free  set 
       $A\subs {\kappa}$ with $|A|={\lambda}$, and  
\item   $\freeUnion{{{\kappa}}}{{\mu}}{{\lambda}}$ iff for each
set-mapping $F:{\kappa}\to {[{\kappa}]}^{<{\mu}}$, the set ${\kappa}$ 
can be partitioned into ${\lambda}$-many  $F$-free  sets.  
       \end{enumerate} 
    \end{definition}

    Hajnal \cite{Ha60} proved $\freeset{{\kappa}}{{\omega}}{{\kappa}}$ and 
    Fodor \cite{Fo52} proved $\freeUnion{{\kappa}}{{\omega}}{{\omega}}$ 
     for each uncountable  cardinal ${\kappa}\in On$ in ZFC.
 
     \begin{theorem}[ZF]\label{tm:free-union-elementary}
        $\freeset{{\kappa}}{\omega}{\omega}$ and 
        $\freeUnion{{\kappa}}{\omega}{\omega}$ hold for each uncountable cardinal ${\kappa}\in On $.
        \end{theorem}
    
We are able to  prove $\freeUnion{{\kappa}}{\omega}{\omega}$
in ZF using a purely combinatorial argument, 
but to derive $\freeset{{\kappa}}{{\omega}}{{\kappa}}$ from ZF 
we should use  the absoluteness method of Karagila.

    \begin{definition}[\cite{Mi37},\cite{Ko84}]
        A family $\mc A$ is {\em $n$-almost disjoint} iff $|A\cap B|<n$ for each $\{A,B\}\in {[\mc A]}^{2}$.   
The family $\mc A$ has {\em property B} iff there is a set $X$ such that 
$X\cap A\ne \empt\ne A\setm X$ for each $A\in \mc A$.
        We say that  $\mc A$ is {\em essentially disjoint} iff there is a 
        function $f$ with $\dom(f)=\mc A$ and $f(A)\in [A]^{<{\omega}}$
        such that $\{A\setm f(A):A\in \mc A\}$ is a family of pairwise disjoint sets.    
        
        Write $M({\kappa},{\omega},n)\to B$   iff every $n$-almost disjoint 
        family $\mc A\subs {[{\kappa}]}^{{\omega}}$ has property $B$, and 
       write $M({\kappa},{\omega},n)\to ED$ iff every $n$-almost disjoint 
       family $\mc A\subs {[{\kappa}]}^{{\omega}}$ is essentially disjoint.
       \end{definition}

        For each infinite cardinal ${\kappa}$ and natural number $n$
          \begin{enumerate}[(1)]
\item $M({\kappa},{\omega},n)\to B$ was proved   by Miller \cite{Mi37}, and 
\item $M({\kappa},{\omega},n)\to ED$ was proved   by Komjáth \cite{Ko84} 
\end{enumerate}    
in ZFC. Let us observe that Komjáth's result is clearly stronger because 
if a family $\mc A\subs {[{\kappa}]}^{{\omega}}$ is essentially disjoint, then 
it has property B, as well.

In Theorem \ref{tm:zf-model}(\ref{en:B}) we prove 
$M({\kappa},{\omega},n)\to B$ in ZF. On the other hand, 
in section \ref{sc:ind} we obtain the following corollary from theorem \ref{tm:model1}:
\begin{corollary}\label{cr:notED}
ZF + $cf({\omega}_1)={\omega}_1$ does not imply  $M({\omega}_1,{\omega},2)\to ED$.
\end{corollary}

However,   some  strengthening of the assumption of $cf({\omega}_1)={\omega}_1$
will be  enough to prove  $M({\omega}_1,{\omega},n)\to ED$ in ZF.
First we need to recall a definition.

        \begin{definition}[Litman,\cite{Li76}]
            A function $f$  is a {\em uniform denumeration on ${\omega}_1$}
            iff $\dom(f)={\omega}_1$ and $f({\alpha})$ is a function from ${\omega}$ onto ${\alpha}$
            for each $1\le {\alpha}<{\omega}_1$. 
             Let {\em UD$({\omega}_1)$} be the assertion that there exists 
             a uniform denumeration on ${\omega}_1$. 
            \end{definition}

        Litman \cite[Lemma 2.8]{Li76} proved that  
        UD$({\omega}_1)$ 
        implies $cf({\omega}_1)={\omega}_1$. The following theorem improves his result. 
 
        \begin{theorem}[ZF]\label{tm:omega1}
            (1) UD$({\omega}_1)$ implies that  
             $M({\omega}_1,{\omega},n)\to ED$  for each natural number $n\in {\omega}$. 
            \\
            \noindent (2) $M({\omega}_1,{\omega},2)\to ED$ implies $cf({\omega}_1)={\omega}_1$.
            \end{theorem}

In Theorem \ref{tm:ud2ed-gen}
we also show that ZF + UD$({\omega}_1)$ implies $M({\kappa},{\omega},n)$ 
for each infinite cardinal ${\kappa}$, but to do so we should create a method using absoluteness 
which makes possible to prove  certain ZFC results using only  ZF + UD$({\omega}_1)$  
(see Theorem \ref{tm:inner-oo}).

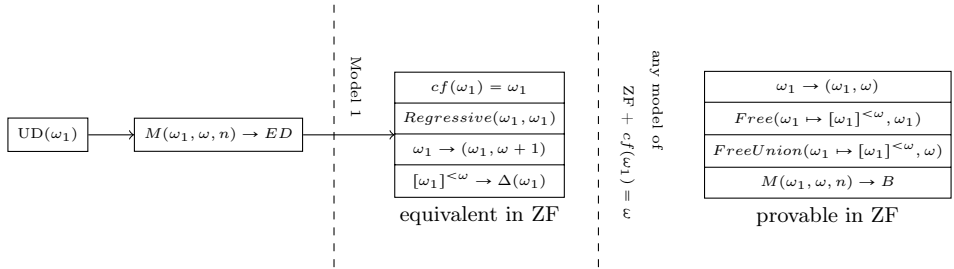
\begin{figure}[ht!]
    \caption{\small Relationship  between different combinatorial properties of ${\omega}_1$}\label{fig:summ}
    \begin{tikzpicture}
    [scale=0.9,every node/.style={scale=0.9},pro/.style={rectangle, inner sep=4pt,minimum size=4mm,draw=black},
    prosplit/.style=  {inner sep=4pt,minimum size=4mm,draw=black,rectangle
    split,rectangle split parts=2,
    rectangle split part fill={white, black!10}}]
    
    \matrix[row sep=5mm,column sep=3mm]{
&&&\node (t1) {};&&\node (t11) {};&&\\
\snode {a1}{UD($\omega_1$)};&&\snode {b1}{$M({\omega}_1,{\omega},n)\to ED$};&
\node[label=above right:\rotatebox{-90}{\tiny Model 1}] (t32) {};&
\node[draw,rectangle,label=below:\small equivalent in ZF,
rectangle split,rectangle split parts=4] (c1) {
\nodepart{one}      \tiny $cf({\omega}_1)={\omega}_1$
\nodepart{two}    \tiny $Regressive({\omega}_1,{\omega}_1)$
\nodepart{three}    \tiny ${\omega}_1\to({\omega}_1,{\omega}+1)$
\nodepart{four}    \tiny $\dsystem{{\omega}_1}{{\omega}}{{\omega}_1}$}
;&\node[label=right:\rotatebox{-90}{\minibox{\tiny any model of \\\tiny \qquad ZF + $cf({\omega}_1)={\omega}$}}] (t32) {};&
\node[draw,rectangle,label=below:\small provable in ZF,
rectangle split,rectangle split parts=4] (d1) {
\nodepart{three}      \tiny $\freeUnion{{\omega}_1}{{\omega}}{{\omega}}$
\nodepart{two}      \tiny $\freeset{{\omega}_1}{{\omega}}{{\omega}_1}$
\nodepart{one}    \tiny ${\omega}_1\to ({\omega}_1,{\omega})$
\nodepart{four}    \tiny $M({\omega}_1,{\omega},n)\to B$
};\\
&&&\node (t2) {};&
&\node (t22) {};&
&\\
}; 
\path[->]  
(a1) edge (b1)
(b1) edge (c1)
;
\path[-, dashed]  (t1) edge (t2);  
\path[-, dashed]  (t11) edge (t22);  
    \end{tikzpicture}

\end{figure}

Figure \ref{fig:summ} summarizes our ZF results concerning combinatorial properties of ${\omega}_1$. 

\noindent {\bf Folklore Facts (ZF).}
\begin{enumerate}[(F1)]
\item \label{f1} If ${\alpha}$ is an infinite ordinal, then 
$|{\alpha}^{<{\omega}}|=|{\alpha}|$.
\end{enumerate}

A non-empty family ${\mathcal {F}}$ of sets  is of {\em finite character} provided 
$ {[F]}^{<{\omega}}\subs \mc F $ for each $F\in \mc F$, 
and ${[Y]}^{<{\omega}}\subs \mc F$ implies $Y\in \mc F$ for each set $Y$.

\begin{enumerate}[(F1)]\addtocounter{enumi}{1}
\item \label{f3} If $A$ is a set of ordinals 
and $\mc F\subs \mc P(A)$ is of finite character, then 
there is an operation $\Gamma_{\mc F}$ such that 
given any set $B\subs  A$,  $\Gamma_{\mc F}(B)$ is a $\subs$-maximal subset of 
$\mc F\cap \mc P(B)$. (Let $\Gamma_{\mc F}(B)=C$ iff $C\in \mc F$ and 
$(C\cap {\alpha})\cup\{{\alpha}\}\notin \mc F$
for each ${\alpha}\in B\setm C$.) 

\end{enumerate}

\section{Equivalence and weakening}\label{sec:elementary}
\begin{definition}
    Assume that  ${\kappa},{\lambda}\in On$  are infinite cardinals.  Write  
    $|\mbb T_{<{\kappa}}|\le {\lambda}$
    iff for each tree   $T\subs On^{<{\omega}}$ if  
    \begin{displaymath}
    \forall s\in T\ |\{{\alpha}\in On: s^\frown {\alpha}\in T\}|<{\kappa}
    \end{displaymath}
    then   $|T|\le {\lambda}.$
    \end{definition}

    In ZFC we have  
    $cf({\kappa}^+)={\kappa}^+$ and so 
    $|\mbb T_{<{\kappa}^+}|\le {\kappa}$ for each infinite cardinal ${\kappa}$.

\begin{theorem}[ZF]\label{tm:equivalence-in-ZF}
    For each infinite cardinal ${\kappa}\in On$ the following statements are equivalent:
    \begin{enumerate}[(i)]
        \item \label{en:cof-zf-eq} $cf({\kappa}^+)= {\kappa}^+$,  
        \item \label{en:tree-zf-eq}     $|\mbb T_{<{{\kappa}^+}}|\le {\kappa}$, 
        \item \label{en:reg-zf-eq} $Regressive({\kappa}^+,{\kappa}^+)$, 
        \item \label{en:delta-zf-eq} $\dsystem{{\kappa}^+}{{\omega}}{{\kappa}^+}$. 
    \end{enumerate}
    Moreover, $cf({\kappa}^+)= {\kappa}^+$ implies that 
    \begin{enumerate}[(i)]\addtocounter{enumi}{4}
        \item \label{en:edm+-zf-eq} ${\kappa}^+\to ({\kappa}^+,{\omega}+1)$,
    \end{enumerate}
    and (\ref{en:edm+-zf-eq}) implies that     
    \begin{enumerate}[(i)]\addtocounter{enumi}{5}
        \item \label{en:cofom-zf-eq} $cf({\kappa}^+)>{\omega}$.
    \end{enumerate}
    \end{theorem}

    \begin{theorem}[ZF]\label{tm:elementary-in-ZF}
        For each infinite cardinal ${\kappa}\in On$ we have 
        \begin{enumerate}[(1)]
        \item \label{en:cof-zf}$ {\kappa}^+\to ({\kappa})^1_{\kappa}$, 
        \item \label{en:tree-zf}     $|\mbb T_{<{{\kappa}}}|\le {\kappa}$, 
        \item \label{en:reg-zf} $Regressive({\kappa}^+,{\kappa})$, 
        \item \label{en:delta-zf} $\dsystem{{\kappa}^+}{{\omega}}{{\kappa}}$, 
        \item \label{en:edm+-zf} ${\kappa}^+\to ({\kappa},{\omega}+1)$.
        \end{enumerate}
        \end{theorem}

First we need a lemma.

\begin{lemma}\label{lm:4in1}\label{lm:cf2trww}
    If ${\kappa},{\lambda}\in On $  are infinite cardinals, ${\lambda}\ge {\kappa}$, 
    and 
    \begin{enumerate}[(a)]
      \item   ${\kappa}^+\to ({\lambda})^1_{\kappa}$,
    \end{enumerate} 
    then 
    \begin{enumerate}[(a)]\addtocounter{enumi}{1}
        \item     $|\mbb T_{<{\lambda}}|\le {\kappa}$, 
        \item     $\dsystem{{\kappa}^+}{{\omega}}{{\lambda}}$,
\item $Regressive({\kappa}^+,{\lambda})$,

\item ${\kappa}^+\to ({\lambda},{\omega}+1)$.
    \end{enumerate}
\end{lemma}

\begin{proof}
$(a)\to(b)$

Let     
 $T\subs {[On]}^{<{\omega}}$ be a  tree such that 
$|\{{\alpha}\in On:s^\frown {\alpha}\in T\}|<{\lambda}$ 
for each $s\in T$.

We should distinguish two cases: since ${\kappa}^+\to ({\lambda})^1_{\kappa}$
and ${\lambda}\ge {\kappa}$, we have either ${\lambda}={\kappa}$ or ${\lambda}={\kappa}^+$.

Assume first that  ${\lambda}={\kappa}^+$.  

For $n\in {\omega}$ write  $T_n=\{s\in T: |s|={\omega}\}$.
By induction on $n$, applying the assumption (a) 
we obtain that $|T_n|\le {\kappa}$.
Thus,  $T$ is the union of ${\omega}$-many sets of size at most 
${\kappa}$.  So $|T|\le {\kappa}$  applying 
(a) again.

Consider now that case when ${\lambda}={\kappa}$.
Then we can define an injective function 
${\varphi}:T\to {\kappa}^{<{\omega}}$ as follows:
\begin{displaymath}
{\varphi}(\<{\alpha}_i:i<n\>)=\<{\beta}_i:i<n\>,
\end{displaymath} 
where 
\begin{displaymath}
{\beta}_j=tp\{{\gamma}:{\gamma}<{\alpha}_j\land 
\<{\alpha}_i:i<j\>^\frown {\gamma}\in T\}.
\end{displaymath}

Since $|{\kappa}^{<{\omega}}|={\kappa}$, we proved 
$|T|\le {\kappa}$.

\smallskip

\noindent $(b)\to (c)$

Assume on the contrary that $\mc A\subs {[On]}^{<{\omega}}$ is a family with 
$|\mc A|={\kappa}^+$  
such that  if $\mc B\subs \mc A$ is a $\Delta$-system, then $|\mc B|< {\lambda}$.

Let $X=\bigcup A$. Since $X$ is a set of ordinals, it has a natural well-ordering,
and so ${[X]}^{<{\omega}}$  also has a natural well-ordering.

Let 
\begin{displaymath}
\mbb F=\{\mc C\subs {[X]}^{<{\omega}}:   \text{the elements of $\mc C$ are pairwise disjoint.}\}
\end{displaymath}
Clearly, $\mbb F$ is of finite character. Thus, by (F\ref{f3}) there is an operation $\Gamma$ such that $\Gamma(\mc D)$ is a $\subs$-maximal 
element of $\mbb F\cap \mc P(\mc D)$ for each $\mc D\subs {[X]}^{<{\omega}}$, i.e. $\Gamma(\mc D)$ is a maximal disjoint subfamily of 
$\mc D$.

Define a tree $T\subs X^{<{\omega}}$ as follows.

We construct the levels $\<T_n:n<{\omega}\>$  of the tree by recursion. 

Let $T_0=\{\empt\}$.

Assume that we have defined $T_n\subs X^n$.

For each $t\in T_n$ do the following. Let 
$$\mc A_t=\{A\setm \ran t:\ran t\subs A\in \mc A\}\setm \{\empt\},$$
and  $\mc B_t=\Gamma(\mc A)$, i.e. $\mc B_t$ is a maximal disjoint subfamily of $\mc A_t$.
Since $\{A\in \mc A: A\setm \ran t\in \mc B_t\}$ is a $\Delta$-system with kernel $\ran t$,
it follows that $|\mc B_t|<{\lambda}$. Since the elements of $\mc B_t$ are finite, $|\bigcup \mc B_t|<{\lambda}$ as well.

To finish the recursive step let 
\begin{displaymath}
T_{n+1}=    
\{t^\frown {\beta}:t\in T_n, {\beta}\in \bigcup \mc B_t\}.
\end{displaymath}

Then $T=\bigcup_{n<{\omega}}T_n$ is a tree and 
\begin{displaymath}
    \forall s\in T\ |\{{\alpha}\in On: s^\frown {\alpha}\in T\}|=|\bigcup \mc B_t|<{\lambda}.
    \end{displaymath}
Thus, $|T|\le  {\kappa}$ by the assumption $|\mbb T_{<{\lambda}}|\le {{\kappa}}$. 
Since $|{\kappa}^{<{\omega}}|={\kappa}$ by (F\ref{f1}), it follows 
$|{[{\kappa}]}^{<{\omega}}|={\kappa}$ and so  
there is $A\in \mc A$ such that $A\not\subset \bigcup \{\ran(t):t\in T\}=\bigcup\{\bigcup \mc B_t:t\in T\}$.

By induction on ${\omega}$ we can 
define a sequence $\<{\beta}_n:n<{\omega}\>$ such that $\<{\beta}_i:i<n\>\in T$ and ${\beta}_i\in A$ as follows.
Assume that $t_n=\<{\beta}_i:i<n\>\in T$ and ${\beta}_i\in A$ for $i<n$.

Then $\empt\ne A\setm \{{\beta}_i:i<n\}\in \mc A_{t_n}$ so ${\beta}_n=\min (A\cap \bigcup \mc B_{t_n})$ is defined by the maximality of $\mc B_{t_n}$.

Thus, $A$ is infinite. Contradiction.

\medskip

\noindent $(b)\to (d)$

Let $f$ be a regressive function on ${\kappa}^+$.
Consider the following tree $T$:
    \begin{multline*}
        T=\big\{s\in ({\kappa}^+)^{<{\omega}}: s=\empt \lor\\ \big (s(0)=0\land   
        \ s(i-1)=f(s(i))\text{ for each } 1\le i<|s|\ \big) \big\}.
    \end{multline*}

Since for each ${\alpha}\in {\kappa}^+$ there is $s\in T$ with 
$s(|s|-1)={\alpha}$,
we have  $|T|={\kappa}^+$.  

Thus, applying the assumption $|\mbb T_{<{\lambda}}|\le {\kappa}$, we obtain that 
there is $s\in T$ such that $B=\{{\beta}\in {\lambda}^+: s^\frown{\beta}\in T \}$
has cardinality at least ${\lambda}$.
If ${\beta}\in B$, then $f({\beta})=s(|s|-1)$ by the construction of $T$.
So taking ${\alpha}=s(|s|-1)$ we have $|f^{-1}\{{\alpha}\}|\ge {\lambda}$.

\medskip

\noindent $(b)\to (e)$
Assume that $f:{[{{\kappa}^+}]}^{2}\to 2$.
    
Define a tree $ T\subs {({{\kappa}^+})}^{<{\omega}}$ as follows:  

By induction on $n$ we define the $n$th level $T_n$ of $T$ as follows:

Let $T_0=\{\empt\}$. 

Assume that $T_{n-1}$  is given. 

For each  $t\in T_{n-1}$ let 
$$A_t=\{{\alpha}\in {{\kappa}^+}: \forall {\zeta}\in ran(t)\  f({\zeta},{\alpha})=1\}.$$
Let $B_t\subs A_t$ be the unique maximal 0 homogeneous set such that 
\begin{displaymath}\tag{$\dag$}\label{dag}
\forall {\alpha}\in  A_t\setm B_t\ (B_t\cap {\alpha})\cup\{{\alpha}\} \text{ is not 
0 homogeneous,}
\end{displaymath}
i.e. $B_t$ was obtained by the  greedy algorithm.

 We can assume that $|B_t|<{\lambda}$, or we proved the statement.

Let $T_{n}=\bigcup \{t^\frown {\beta}:{\beta}\in B_t:t\in T_{n-1}\}$.

Then $T$ is $<{\lambda}$-branching, so $|T|\le {\kappa}$ by the assumption of the Lemma. 
Thus, $|\bigcup\{\ran(t):t\in T\}|\le {\lambda}$, because this set is the union of at most ${\kappa}$ many 
finite sets.

So we can pick      ${\alpha}\in {\kappa}^+\setm \bigcup\{\ran(t):t\in T\}$.

By recursion on ${\omega}$, we can construct a sequence
$\<{\beta}_n:n<{\omega}\>$ such that 
\begin{displaymath}\tag{$\circ_n$}
\forall n\in {\omega}\ \<{\beta}_i:i<n\>\in T_n \land {\alpha}\in A_{\<{\beta}_i:i<n\>}
\land\forall i<n\  f({\beta}_i,{\alpha})=1.
\end{displaymath}
For $n=0$,  $\empt\in T_0$ and   $A_{\empt}={\kappa}$,  so $\circ_0$ holds. 

Assume that we have constructed  $\<{\beta}_i:i<n-1\>$
such that $(\circ_{n-1})$ holds.  Then  ${\alpha}\in A_{\<{\beta}_i:i<n-1\>}\setm 
B_{\<{\beta}_i:i<n-1\>}$, so it follows from (\ref{dag}) that 
\begin{displaymath}
{\beta}_{n-1}=\min\{{\beta}\in  B_{\<{\beta}_i:i<n-1\>}\cap {\alpha}:f({\beta}_{n-1},{\alpha})=1\}
\end{displaymath}
is defined. Then $(\circ_n)$ holds.
So we can carry out the construction. 

Thus, $\{{\beta}_n:n\in {\omega}\}\cup\{{\alpha}\}$ is 1-homogeneous and has order type 
${\omega}+1$ because ${\beta}_n<{\alpha}$ for each $n\in {\omega}$.
\end{proof}

\begin{proof}[Proof of Theorem \ref{tm:elementary-in-ZF}]
  
    According to \cite[Section 30 Problem 27]{KoTo06},  
    the partition relation $ {\kappa}^+\to ({\kappa})^1_{\kappa}$ was proved first by Jech 
    in \cite{Je82} for 
    ${\kappa}={\omega}_1$. To make this note self-contained we include a simple proof.  
    
    Assume that ${\delta}\not\to ({\kappa})^1_{\kappa}$, i.e. 
    there is a function  $f:{\delta}\to {\kappa}$  such that 
    $|f^{-1}\{{\alpha}\}|<{\kappa}$ for each ${\alpha}<{\kappa}$.
    Then  we can construct a one-to-one function  
    ${\varphi}:{\delta}\to {\kappa}\times {\kappa}$ as follows: 
        $$
        {\varphi}({\alpha})=\<f({\alpha}),tp\big(\big\{{\gamma}<{\alpha}:f({\gamma})=f({\alpha})\big\}\big)\>.
        $$
    Recalling that $|{\kappa}\times {\kappa}|={\kappa}$ in ZF 
    by (F\ref{f1}) we obtain that $|{\delta}|\le {\kappa}$ and so ${\delta}<{\kappa}^+$.
    Thus, ${\kappa}^+\to ({\kappa})^1_{\kappa}$.

    (\ref{en:cof-zf}) $\to $ (\ref{en:tree-zf}  (\ref{en:reg-zf})(\ref{en:delta-zf}) 
    (\ref{en:edm+-zf}) follow clear from the lemma taking ${\lambda}={\kappa}$.
        \end{proof}

    \begin{proof}[Proof of Theorem \ref{tm:equivalence-in-ZF}]

    \makebox[1pt]{}
\smallskip

        \noindent {\bf Claim}: {\em If $cf({\kappa}^+)={\kappa}^+$, then 
        ${\kappa}^+\to ({\kappa}^+)_{\kappa}$. 
    }
    \begin{proof}[Proof of the Claim]
    Assume that  
        $\mc A=\{A_n:n\in {\kappa}\}\subs [{{\kappa}^+}]^{{\kappa}}$. 
If $A\in \mc A$, then $A$ is not cofinal in ${\kappa}^+$, so $sup(A)<{\kappa}^+$   
        The set    $\{\sup(A):A\in \mc A\}\in {[{\kappa}^+]}^{{\kappa}}$, and so 
        ${\beta}=\sup B<{\kappa}^+$. Thus, $|\bigcup \mc A|\le |{\delta}|<{\kappa}^+$.
   Thus, ${\kappa}^+\to ({\kappa}^+)_{\kappa}$ holds. 
    \end{proof}

(\ref{en:cof-zf-eq}) $\to $ $(\ref{en:tree-zf-eq})\land (\ref{en:reg-zf-eq})\land (\ref{en:delta-zf-eq})
\land (\ref{en:edm+-zf-eq})$ follows  from the    lemma \ref{lm:4in1} because 
it can be applied for  
${\kappa}={\lambda}$
by the  Claim.

If (\ref{en:cof-zf-eq}) fails, we will give counterexamples
for (\ref{en:tree-zf-eq}),  (\ref{en:reg-zf-eq}) and (\ref{en:delta-zf-eq}).
Assume that $\{{\alpha}_{\zeta}:{\zeta}<{\kappa}\}$ is cofinal in ${\kappa}^+$.
We can assume that ${\alpha}_0={\kappa}$

Let 
$T=\{\<{\zeta},{\eta}\>:{\zeta}\in {\kappa}, {\eta}<{\alpha}_{\zeta}\}$.
Then $T$ witnesses that $|\mbb T_{<{\kappa}^+}|\le {\kappa}$ fails. 

Next 
define a regressive function $g$ on ${\kappa}^+$ as follows:
$g(0)=0$ and for ${\alpha}>0$,
\begin{displaymath}
g({\alpha})=\min\{{\xi}<{\kappa}:{\alpha}\le {\alpha}_{\zeta}\}.
\end{displaymath}
Then $g^{-1}\{{\zeta}\}\subs {\alpha}_{\zeta}$ for each ${\zeta}<{\kappa}$.
Thus,  $g$ shows that $Regressive({\kappa}^+, {\kappa}^+)$ fails.

Finally,  let  
         $$\mc A=\big\{\{{\zeta},{\alpha}\}:{\zeta}<{\omega},{\kappa}\le {\alpha}<{\alpha}_{\zeta}\big\}\subs {[{\omega}_1]}^{2}.$$
         Then $\mc A$ has cardinality ${\kappa}^+$, but it does not contain a 
         $\Delta$-system of size ${\kappa}^+$, 
         so  $\dsystem{{\kappa}^+}{{\omega}}{{\kappa}^+}$ fails. 
    
So we proved 
$(\ref{en:tree-zf-eq})\lor  (\ref{en:reg-zf-eq}) \lor  (\ref{en:delta-zf-eq})\to 
(\ref{en:cof-zf-eq})$. 

What remained is to show that (\ref{en:edm+-zf-eq}) implies (\ref{en:cofom-zf-eq}). 
Assume that $cf({\kappa}^+)={\omega}$ and let $\<{\alpha}_n:n<{\omega}\>$
is a strictly increasing cofinal sequence in ${\kappa}$    
Define $f:{[{\omega}_1]}^{2}\to 2$ as follows: for ${\zeta}<{\xi}<{\kappa}^+$ let
$f({\zeta},{\xi})=1$ iff ${\zeta}<{\alpha}_n\le {\xi}$ for some $n<{\omega}$.
Then $f$ proves that ${\kappa}^+\not\to ({\kappa}^+,{\omega}+1)$. 
\end{proof}

    \section{More results with purely combinatorial proofs} \label{sc:more}

\begin{proof}[Proof of Theorem \ref{tm:free-union-elementary}.  Part 1] 
We show that   $\freeUnion{{\kappa}}{\omega}{\omega}$ hold for each uncountable cardinal ${\kappa}\in On $
using a purely combinatorial argument. 

Assume  that $F:{\kappa}\to {[{\kappa}]}^{<{\omega}}$.

    Consider the sequence  $\<B_{\alpha}:{\alpha}<{\kappa}\>$ such that 
    ${\alpha}\in B_{\alpha}$ and $B_{\alpha}$ is the minimal $F$-closed subset of ${\kappa}$
    which contains ${\alpha}$, i.e. ${\alpha}\in B_{\alpha}$ and 
    ${\zeta}\in B_{\alpha}$ implies $F({\zeta})\subs B_{\alpha}$. 
    Since $B_{\alpha}$ can be obtained by recursion  
    as an increasing union of finite sets, $|B_{\alpha}|\le {\omega}$.
    Moreover, we can obtain  
a function $e$ with $\dom(e)={\kappa}$
    such that $e({\alpha})$ maps ${\omega}$ onto $B_{\alpha}$. 
    
    For ${\alpha}<{\kappa}$ write $B_{<{\alpha}}=\bigcup_{{\zeta}<{\alpha}}B_{{\zeta}}$
    and $A_{\alpha}=B_{\alpha}\setm B_{<{\alpha}}$.
    
    We will define a function $g:{\kappa}\to {\omega}$ such that 
    $g^{-1}\{n\}$ is $F$-free
     for each 
     $n<{\omega}$.

    By transfinite recursion we define functions $g_{\alpha}$ such that $\dom(g_{\alpha})=A_{\alpha}$
    and $g=\bigcup_{{\alpha}<{\kappa}}g_{\alpha}$ satisfies the requirements.
    
    Assume that we have defined the functions $\<g_{\zeta}:{\zeta}<{\alpha}\>$.
    Write  $g_{<{\alpha}}=\bigcup_{{\zeta}<{\alpha}}g_{\zeta}.$
    Then $g_{<{\alpha}}:B_{<\alpha}\to {\omega}$.
    
    If $A_{\alpha}=\empt$, then let $g_{\alpha}=\empt$.
    
    If $A_{\alpha}\ne \empt $, then 
    $e({\alpha})$ maps ${\omega}$ onto $B_{\alpha}$,
    so we can construct  $h({\alpha})$ mapping ${\omega}$ onto $A_{\alpha}$.

    By induction on $n\in {\omega}$
    define $g'_{\alpha}:{\omega}\to {\omega}$ as follows: 
    \begin{multline*}\tag{$\ddag$}
    g'_{\alpha}(n)=\\\min\big({\omega}\setm 
    (\{g'_{\alpha}(m):m<n\}\cup \{g_{<{\alpha}}({\zeta}):{\zeta}\in F(h_{\alpha}(n))\cap B_{<{\alpha}}\})\big).
    \end{multline*}
    The definition  of $g'_{\alpha}(n)$ is meaningful, because  
    $\{g'_{\alpha}(m):m<n\}\cup \{g_{<{\alpha}}({\zeta}):{\zeta}\in F(h_{\alpha}(n))\cap B_{<{\alpha}}\}$ 
    is a finite set.

    Then, for ${\eta}\in A_{\alpha}$ let
    \begin{displaymath}
    g_{\alpha}({\eta})=g'_{\alpha}(n) \text{ where }
    n=\min\{m\in {\omega}:h_{\alpha}(m)={\eta}\}.
    \end{displaymath}
    The definition  $g_{\alpha}({\eta})$ is meaningful because $h_{\alpha}$ is onto. 
    
    By the construction, $g^{-1}\{k\}$ is $F$-free for each $k<{\omega}$.
    Indeed, assume that $\{{\xi},{\eta}\}\in {[{\kappa}]}^{2}$.
    Pick ${\alpha},{\beta}\in {\kappa}$ such that ${\xi}\in A_{\alpha}$ and ${\eta}\in A_{\beta}$.
    
    If ${\alpha}={\beta}$, then $g({\xi})=g_{\alpha}({\xi})=g'_{\alpha}(n)$
    and $g({\eta})=g_{\alpha}({\eta})=g'_{\alpha}(m)$ 
    for some $n,m<{\omega}$ with 
    $h_{\alpha}(m)={\eta}$ and $h_{\alpha}(n)={\xi}$. Thus, $n\ne m$. Since 
    $g'_{\alpha}$ is injective, 
    $g({\eta})=g_{\alpha}({\eta})\ne g_{\alpha}({\xi})=g({\xi})$.
    
    If ${\alpha}\ne {\beta}$, then we can assume next that ${\beta}< {\alpha}$. Since $B_{\beta}$ is $F$-closed,
    ${\xi}\notin F({\eta})$. 
    Assume  ${\eta}\in F({\xi})$. Then $g_{\alpha}({\xi})=g'_{\alpha}(n)$ for some 
    $n\in {\omega}$ with ${\xi}=h_{\alpha}(n)$.  But then  
    ${\eta}\in  F(h_{\alpha}(n))\cap B_{<{\alpha}}$
    and so $g_{\alpha}'(n)\ne g_{<{\alpha}}({\eta})=g({\eta})$ by $(\ddag)$.  
\end{proof}

    \begin{proof}[Proof of Theorem \ref{tm:omega1}]
    $(1)$. 
    Assume that $\mc A\subs [{\omega}_1]^{\omega}$ is an $n$-almost disjoint family.
    
    Since $|{[{\omega}_1]}^{n}|={\omega}_1$ we have $|\mc A|\le {\omega}_1$.
    So we can assume that $\bigcup \mc A={\omega}_1$ and 
    $\mc A=\{A_{\alpha}:{\alpha}<{\omega}_1\}$.

    Let
    \begin{displaymath}
    C=\{{\gamma}<{\omega}_1:{\gamma}=\bigcup_{{\beta}<{\gamma}}A_{\gamma}\land \forall {\delta}\in {\omega}_1\setm {\gamma}
    \ |A_{\delta}\cap {\gamma}|<n\}.
    \end{displaymath}
    Then $C$ is closed and unbounded  in ${\omega}_1$ because 
    $cf({\omega}_1)={\omega}_1$. Let $\{{\gamma}_{\zeta}:{\zeta}<{\omega}_1\}$
    be the unique increasing enumeration of $C$. Clearly ${\gamma}_0=0$.
    
    Let $f$ be a uniform enumeration on ${\omega}_1$.
    Using $f$ as a parameter, we can construct a bijection $h:{\omega}_1\to {\omega}_1$
    such that $\{h({\omega}{\xi}+n):n<{\omega}\}=[{\gamma}_{\xi},{\gamma}_{{\xi}+1})$.

    
    Define $F:{\omega}_1\to {[{\omega}_1]}^{<{\omega}}$ as follows:
    \begin{displaymath}
    F(h({\omega}{\xi}+n))=A({\omega}{\xi}+n)\cap \big({\gamma}_{\xi}\cup \bigcup 
    \{A(h({\omega}{\xi}+m):m<n\}\big).
    \end{displaymath}
    It is easy to see that the sets $\{A_{{\alpha}}\setm F({\alpha}):{\alpha}<{\omega}_1\}$
    are pairwise disjoint.
        So $\mc A$ is essentially disjoint.

    \medskip\noindent (2)
    
    We prove the contrapositive. 
    
    Assume that $\{{\alpha}_n:n<{\omega}\}$ is cofinal in ${\omega}_1$.
    We can assume that ${\alpha}_0=0$ and  ${\alpha}_i+{\alpha}_i<{\alpha}_{i+1}$
    for each $i<{\omega}$. 
    
    Define $\mc A\subs {[{\omega}_1]}^{{\omega}}$ as follows. 
    For $n<{\omega}$ let 
    \begin{displaymath}
    I_n={\alpha}_{n+1}\setm {\alpha}_n,\end{displaymath}
    and for $1\le n <{\omega}$ and ${\alpha}_{n-1}\le {\zeta}<{\alpha}_n$ let
    \begin{displaymath}
    A(n,{\zeta})=\{{\alpha}_{n+i}+{\zeta}:i\in {\omega}\}.
    \end{displaymath}
    Put 
    \begin{displaymath}
    \mc A=\{I_n:n<{\omega}\}\cup\{A(n,{\zeta}):1\le n<{\omega},{\alpha}_{n-1}\le {\zeta}<{\alpha}_n\}.
    \end{displaymath}
    The family $\mc A$ is clearly $2$-almost disjoint.
    Indeed, the $I_n$ are pairwise disjoint and clearly  $|I_m\cap A(n,{\zeta})|\le 1$.
    If ${\rho}\in I(n,{\zeta})\cap I(m,{\eta})$, then 
    ${\rho}={\alpha}_{n+i}+{\zeta}={\alpha}_{m+j}+{\eta}$ for some $i,j\in {\omega}$.
    Since ${\alpha}_{n+i}+{\zeta}<{\alpha}_{n+i}+{\alpha}_{n+i}<{\alpha}_{n+i+1}$, it follows that 
    $n+i=m+j$ and so ${\zeta}={\eta}$. 
    So $I(n,{\zeta})= I(m,{\eta})$.
    Thus, the sets $\{A(n,{\zeta}):1\le n<{\omega},{\alpha}_{n-1}\le {\zeta}<{\alpha}_n\}$
    are pairwise disjoint.
    
    Finally, we show that $\mc A$ is not ED.
    Assume that $f$ is a function with $\dom(f)=\mc A$ and $\ran (f)\subs {[{\omega}_1]}^{<{\omega}}$.
    Let $B=\bigcup\{f(I_n):n<{\omega}\}$. Then $B$ is a countable union of finite sets, so $B$ is countable 
    by Theorem \ref{tm:elementary-in-ZF}(\ref{en:cof-zf}).
    Let $C=\{{\beta}\dotdiv{\alpha}_n:n<{\omega}, {\beta}\in B, {\alpha}_n<{\beta}\}$, where
    ${\beta}\dotdiv{\alpha}_n$ is the unique ordinal ${\zeta}$ with ${\alpha}_n+{\zeta}={\beta}$.
    Then $C$ is also countable because we can map ${\omega}\times {\omega}$ onto $C$. 
    
    Let ${\zeta}\in {\omega}_1\setm C$.  Fix $n$ with ${\alpha}_{n-1}\le {\zeta}<{\alpha}_n$.
    Then $C\cap A(n,{\zeta})=\empt$. Let $m\in {\omega}$ such that 
    ${\alpha}_m+{\zeta}\in A(n,{\zeta})\setm f(A(n,{\zeta}))$. Then 
    $${\alpha}_m+{\zeta}\in (I_m\setm f(I_m))\cap (A(n,{\zeta})\setm f(A(n,{\zeta})).$$
    So $f$ does not witness that $\mc A$ is ED.
    
    Since $f$ was arbitrary, we verified that $\mc A$ is not ED. Thus, we proved  (2) as well. 
    \end{proof}

\begin{definition}\label{df:Hnew}
    Given a set  $X$, let 
    $h(X)$ be the minimal ordinal ${\alpha}>0$ such that 
    there is no function mapping $X$ onto ${\alpha}$.
    \end{definition}

 Let us remark that $h(X)$ is defined, and  $h(X)$ is at most the 
 Hartog's number of $\mc P(X)$.   

    \begin{theorem}[ZF]\label{tm:er}
        For each cardinal ${\kappa}\in On$ and ${\mu}\in On$ there is a cardinal ${\lambda}\in On$
        such that ${\lambda}\to ({\kappa})^2_{\mu}$.
            \end{theorem}
            
        \begin{proof}
    Pick a cardinal ${\sigma}$ with ${\sigma}\to ({\kappa})^1_{\mu}$, and let 
    ${\lambda}=h({\mu}^{<{\sigma}})$.
    
            Assume $c:{[{\lambda}]}^{2}\to {\mu}$.

        Define a partial function $f:{\mu}^{<{\sigma}}\to {\lambda}$ as follows.
        
        Let $f(\empt)=0$.
        
        Let $s\in {\mu}^{<{\sigma}}$. Assume that we have defined $f(s\restriction {\zeta})$ for each ${\zeta}\in \dom(s)$.
        Let  
        \begin{displaymath}
        A_s=\{{\alpha}\in {\lambda}: \forall {\zeta}\in \dom(s)\ c({\alpha},f(s\restriction {\zeta}))=s({\zeta})\}.
        \end{displaymath}
        Let $f(s)=\min  A_s$ if $A_s\ne \empt$. If $A_s=\empt$, then let $f(s)=0$. 
        
        By the choice of ${\lambda}$, we have 
        $\ran f\subsetneq  {\lambda}$, so we can pick 
         ${\alpha}\in {\lambda}\setm \ran f$.
        Then there is $s\in {\mu}^{\sigma}$ such that ${\alpha}\in A_{s\restriction {\zeta}}$ for ${\zeta}<{\sigma}$.
        
        Then $\{f(s\restriction {\zeta}):{\zeta}<{\sigma})\}$ is end-homogeneous.
       Since ${\sigma}\to ({\kappa})^1_{\mu}$, $\{f(s\restriction {\zeta}):{\zeta}<{\sigma})\}$
       contains a $c$-homogeneous set  of size ${\kappa}$.
        \end{proof}

\section{Using absoluteness}\label{sc:inner}

In \cite[Example 1]{Ka14}, Karagila demonstrated that the Erdos-Dushnik-Miller Theorem can be proven in ZF. In that paper, specifically 
in \cite[Theorem 3 and 5]{Ka14}, he presented a method for deriving ZF results from ZFC results using absoluteness.
 We believe that the following restatement of his ideas is highly 
 applicable and also reveals the limitations of that approach.

\begin{theorem}\label{tm:inner}
Assume that the formula ${\varphi}(\vec x)$ is downward absolute and the formula 
$\psi(\vec x,\vec y)$ is  
upwards absolute
between transitive ZF models $M\subs N$ with $On^N=On^M$,
and ${\varphi}(a_0,\dots, a_{n-1})$ implies  $a_i\subs L$ for $i<n$, 
If 
\begin{displaymath}
ZFC \vdash \forall \vec x \big ( {\varphi}(\vec x) \to  \exists \vec y 
\psi(\vec x,\vec y) \big ),
\end{displaymath}
then 
\begin{displaymath}
    ZF \vdash \forall \vec x  \big ( {\varphi}(\vec  x) \to 
    \exists \vec y \psi(\vec  x,\vec y) \big ),
    \end{displaymath}
    \end{theorem}

\begin{proof}
Assume that $V\vDash{\varphi}(a_0,\dots, a_{n-1})$. 
Since  $a_0,\dots, a_{n-1}\in L[a_0,\dots, a_{n-1}]$
and $L[a_0,\dots, a_{n-1}]\vDash ZFC$ we have 
 \begin{displaymath}
    L[a_0,\dots, a_{n-1}] \vDash \exists \vec y \psi(a_0,\dots, a_{m-1},\vec y).
\end{displaymath}
 Pick $b_0,\dots, b_m\in L[a_0,\dots, a_{n-1}]$ with 
\begin{displaymath}
    L[a_0,\dots, a_{n-1}] \vDash \psi(a_0,\dots, a_{m-1}, b_0,\dots, b_{m-1}).
\end{displaymath}
Since the formula $\psi(\vec x,\vec y)$ is upwards absolute, we have  
$$V\vDash\psi(a_0,\dots, a_{m-1}, b_0,\dots, b_{m-1}).$$
\end{proof}

\begin{theorem}\label{tm:zf-model}
    For each  infinite cardinal ${\kappa}\in On $ we have 
    \begin{enumerate}[(1)]

    \item ${\kappa}\to ({\kappa},{\omega})$
    (\cite[Example 1]{Ka14})\smallskip
\item if ${\kappa}$ is a regular, then 
${\kappa}\to ({\kappa},{\omega}+1)$, \smallskip
    \item 
    $\freeset{{\kappa}}{{\omega}}{{\kappa}}$, \smallskip 
    \item $\freeUnion{{\kappa}}{{\omega}}{{\omega}}$, \smallskip
    \item\label{en:B}  $M({\kappa},{\omega},n)\to B$. \smallskip
\end{enumerate}

\end{theorem}

\begin{proof}
(2) Let ${\varphi}({\kappa},f)$ be the formula 
\begin{displaymath}
\text{${\kappa}\in On$ is an infinite regular cardinal and $f:{[{\kappa}]}^{2}\to 2$},
\end{displaymath}
and let ${\psi}({\kappa},f,A,g)$ be the formula 
\begin{multline*}
\text{$f:{[{\kappa}]}^{2}\to 2$, $A\subs {\kappa}$ and (either $g:{\omega}\to A$ is a bijection and 
$f''[A]^2=\{1\}$) or }\\
\text{ $g:{\kappa}\to A$ is a bijection and 
$f''[A]^2=\{0\}$)},
\end{multline*}
We can apply Theorem \ref{tm:inner}.

\noindent (3) and (4) can be proved similarly. 

\noindent (5)
We need some preparation because the statement 
``$A$ is countable'' is not absolute. However,  the following statement clearly 
implies $M({\kappa},{\omega},n)\to B$
\begin{enumerate}[(*)]
\item if $\mc A\subs \mc P({\kappa})$ is $n$-almost disjoint and $tp(A)={\omega}$ for each 
$A\in \mc A$, then $\mc A$ has property $B$,
\end{enumerate}
and the formula ``$A\subs On $ has order type ${\omega}$'' is absolute. 
So let ${\varphi}({\kappa},\mc A)$ be the formula 
\begin{multline*}
\text{${\kappa}\in On$ is an infinite cardinal, $\mc A\subs \mc P({\kappa})$
is $n$-almost disjoint,}\\
\text{and $tp(A)={\omega}$ for each $A\in \mc A$},
\end{multline*}
and let ${\psi}({\kappa},\mc A,g)$ be the formula 
\begin{displaymath}
\text{$\mc A\subs \mc P({\kappa})$, $g:{\kappa}\to 2$,  and 
$g''A=\{0,1\}$ for each $A\in \mc A$.  }
\end{displaymath}
\end{proof}

\begin{problem}\label{cr:elementary}
    Find elementary (combinatorial) proofs for Theorem \ref{tm:zf-model}(1)--(5).
    \end{problem}

\noindent
{\bf Remark:} 
We have seen that $M({\omega}_1,{\omega},2)\to ED$ implies $cf({\omega}_1)={\omega}_1$.
We will see that 
$M({\omega}_1,{\omega},2)\to ED$ is not provable even from  ZF + $cf({\omega}_1)={\omega}_1$. 
However, we can get the following result:

\begin{theorem}[ZF]\label{tm:ud2ed-gen}
If $UD({\omega}_1)$ holds, ${\kappa}$ is an arbitrary cardinal, and $n\in {\omega}$,
then every $n$-almost disjoint family $\mc A\subs {[{\kappa}]}^{{\omega}}$ is essentially 
disjoint.  
\end{theorem}

To prove this result  we need  the following 
 theorem,  which is based on the ideas of Karagila's method: it makes possible to
obtain results from ZF + UD(${\omega}_1$) using absoluteness.

\begin{theorem}\label{tm:inner-oo}
    Assume that the formula ${\varphi}(\vec x)$ is downwards absolute  and $\psi(\vec x,\vec y)$ 
    is upwards absolute 
        between transitive ZF models $M\subs N$ with $On^N=On^M$ and ${\omega}_1^N={\omega}_1^M$,
        moreover 
    ${\varphi}(a_0,\dots, a_{n-1})$ implies  $a_i\subs L$ for $i<n$. 
    If 
    \begin{displaymath}
    ZFC \vdash \forall \vec x \big ( {\varphi}(\vec x) \to  \exists \vec y 
    \psi(\vec x,\vec y) \big ),
    \end{displaymath}
    then 
    \begin{displaymath}
        ZF + UD(\omega_1)
        \vdash \forall \vec x  \big ( {\varphi}(\vec  x) \to 
        \exists \vec y \psi(\vec  x,\vec y) \big ).
        \end{displaymath}
        \end{theorem}
    First, we need a technical lemma. 
        \begin{lemma}\label{lm:un-eq}
            UD$({\omega}_1)$ 
            iff 
            there is a set $p$ such that 
            \begin{displaymath}
            {{\omega}_1}^{L[p]}={\omega}_1.
            \end{displaymath}
            \end{lemma}
            
            \begin{proof}
            If $f$ is a uniform denumeration on ${\omega}_1$, then let
            $$p=\{\<{\alpha},n,f({\alpha})(n)\>:0<{\alpha}<{\omega}_1, n\in {\omega}_1\}.$$
Then $p\subs L$, and so  $f\in L[p]$. Thus, every ${\alpha}<{\omega}_1$ is countable in $L[p]$,
and so   ${{\omega}_1}^{L[p]}={\omega}_1$.
            
            If ${{\omega}_1}^{L[p]}={\omega}_1$, then let $f$ be a uniform  denumeration of
            ${{\omega}_1}^{L[p]}$ in $L$. Since  ${{\omega}_1}^{L[p]}={\omega}_1$, 
            $f$ is a  uniform  denumeration of
            ${\omega}_1$ in $V$.
            \end{proof}

        \begin{proof}[Proof of  Theorem \ref{tm:inner-oo}]
            Assume $V\vDash$ ZF + UD($\omega_1$). By Lemma \ref{lm:un-eq}
            there is a set $p$ such that  
            \begin{displaymath}
            {{\omega}_1}^{L[p]}={\omega}_1.
            \end{displaymath}

            Assume that $V\vDash{\varphi}(a_0,\dots, a_{n-1})$. 
            Since  $a_0,\dots, a_{n-1}\in L[a_0,\dots, a_{n-1}]\subs  L[a_0,\dots, a_{n-1},p]$
            and $L[a_0,\dots, a_{n-1},p]\vDash ZFC$ and ${\omega}_1^{L[a_0,\dots, a_{n-1},p]}={\omega}_1$, 
            we have 
            \begin{displaymath}
                L[a_0,\dots, a_{n-1},p] \vDash {\varphi}(a_0,\dots, a_{m-1}),
            \end{displaymath}
and so 
            \begin{displaymath}
                L[a_0,\dots, a_{n-1},p] \vDash \exists \vec y \psi(a_0,\dots, a_{m-1},\vec y).
            \end{displaymath}
             Pick $b_0,\dots, b_m\in L[a_0,\dots, a_{n-1},p]$ with 
            \begin{displaymath}
                L[a_0,\dots, a_{n-1},p] \vDash \psi(a_0,\dots, a_{m-1}, b_0,\dots, b_{m-1}).
            \end{displaymath}
            Since the formula $\psi(\vec x,\vec y)$ is  upwards absolute, we have  
            $$V\vDash\psi(a_0,\dots, a_{m-1}, b_0,\dots, b_{m-1}).$$
            \end{proof}

            \begin{proof}[Proof of Theorem \ref{tm:ud2ed-gen}]
                Let ${\varphi}({\kappa},\mc A)$ be the formula 
                \begin{displaymath}
                \text{${\kappa}\in On$ is an infinite cardinal, 
                $\mc A\subs {[{\kappa}]}^{{\omega}}$
                is $n$-almost disjoint},
                \end{displaymath}
                and let ${\psi}({\kappa},\mc A,g)$ be the formula 
                \begin{displaymath}
                \text{$\mc A\subs {[{\kappa}]}^{{\omega}}$, $g:{\kappa}\to {[{\kappa}]}^{<{\omega}}$,  and 
                $\{A\setm g(A):A\in \mc A\}$ is a disjoint family.  }
                \end{displaymath}
            The formula ``{\em $A\subs On$ is countable}'' is absolute 
            between transitive ZF models $M\subs N$ with $On^N=On^M$ and ${\omega}_1^N={\omega}_1^M$,
            so we can apply Theorem \ref{tm:inner-oo}.
            \end{proof}

\section{Independence results in ZF}\label{sc:ind}

In this section we investigate the implication between  
certain combinatorial statements concerning ${\omega}_1$.

\begin{theorem}\label{tm:model1}
The following two statements are equiconsistent:
\begin{enumerate}[(1)]
\item ZFC+ there is an inaccessible cardinal
\item 
ZF + 
\begin{enumerate}[(i)]
    \smallskip\item $cf({\omega}_1)={\omega}_1$,
\smallskip\item 
${\omega}_1\to ({\nu})^n_{\omega}$
for each ${\nu}<{\omega}_1$ and  $n<{\omega}$,
    \smallskip\item  $M({\omega}_1,{\omega},2)\to ED$ fails,
\smallskip\item there is no universal denumeration of ${\omega}_1$.
\end{enumerate}
\end{enumerate}    
Assuming the existence of a weakly compact cardinal, 
you can get the consistency of ZF + (i), (ii), (iii), (iv) +
\begin{enumerate}[(v)]
    \item ${\omega}_1\to ({\omega}_1)^n_{\omega}$ for each $n<{\omega}$.
\end{enumerate}
    \end{theorem}

In the next  proof we use freely the terminology of 
\cite[Chapter 17:Models in Which AC fails]{Halb2017} concerning symmetric generic extensions.  

If $P$ is a poset, $\underline x$ is a $P$-name and $\mc G\subs P$ is a generic filter over some 
ZFC model $\mc M$, write $\KKGU x$ for the interpretation of $\underline{x}$ in $\mc M[\mc G]$, i.e. 
\begin{displaymath}
    \KKGU x=\{\KKGU y: \exists p\in \mc G\ \<\underline{y},p\>\in \underline{x}\}.
\end{displaymath}

\begin{proof}[Proof of Theorem \ref{tm:model1}]
(1)$\to$ (2). 
    First we construct our model $\mc N_1$.

Assume that $\mc M\models$ ``ZFC + ${\kappa}$ is inaccessible.''
    Let $C=\{{\lambda}<{\kappa}:{\lambda}\text{ is an infinite cardinal}\}$.

Consider the poset $P$ below which collapses every ${\lambda}\in C$ to ${\omega}$:  
\begin{displaymath}
P=\{p\in Fn(C\times {\omega},{\kappa};{\omega}): {\alpha}=p({\lambda},n) \text{ implies }{\alpha}<{\lambda}\},
\end{displaymath}
and let  $p\le q$ iff $p\supset q$.

Next we define a subgroup $G$ of the automorphisms of $P$.
First we define the underlying set of $G$ as follows:
\begin{displaymath}
G={}^CS({\omega}),
\end{displaymath}
i.e. a typical element ${\pi}$ of $G$ is a function which 
assigns a permutation of ${\omega}$ to each infinite cardinal  ${\lambda}$ below ${\kappa}$.
 
If  ${\pi}\in G$ and  $p\in P$ 
define ${\pi}(p)\in P$ as follows:
$$\dom({\pi}(p))=\{ \<{\lambda},{\pi}({\lambda})(n)\>:\<{\lambda},n\>\in \dom(p)\}$$ and     
\begin{displaymath}
    {\pi}(p)({\lambda},{\pi}({\lambda})(n))=p({\lambda},n).
\end{displaymath}
So 
\begin{displaymath}
G\le Aut(P).
\end{displaymath}

For ${\alpha}<{\kappa}$ let
\begin{displaymath}
H_{\alpha}=\{{\pi}\in G: \forall {\lambda}\in C\cap {\alpha}: {\pi}_{\lambda}=id_{\omega}\}.
\end{displaymath}

Let 
\begin{displaymath}
\mc F=\{H\le G: \exists {\alpha}<{\kappa}\ H_{\alpha}\le H \}.
\end{displaymath}
Then 
$\mc F$ is a normal filter. 

Let $HS$ be the class of hereditarily  symmetric names.
Let $\mc G$ be a generic filter in $P$ over $\mc M$, 
and   let 
$$\mc N_1=\{\KKGU x: \underline{x}\in HS\}.$$

We are to show that $\mc N_1$ satisfies the requirements. 

For ${\mu}<{\kappa}$ let 
$P_{\mu}=\{p\in P: \dom(p)\subs {\mu}\times {\omega}\}$.

\begin{lemma}[Key lemma]\label{lm:key}
If ${\alpha},{\beta}\in On$, $f\in \mc N_1$, $f:{\alpha}\to {\beta}$,
then there is ${\mu}<{\kappa}$ such that 
\begin{displaymath}
f\in \mc M[\mc G\cap P_{\mu}].
\end{displaymath}
\end{lemma}
\begin{proof}[Proof of the Key lemma]
Pick $\underline f\in HS$ such that $\KKGU f=f$.
Let ${\mu}<{\kappa}$ such that $st(\underline{f})\ge H_{\mu}$.

\begin{claim}
If $p\Vdash \underline{f}(\check{\zeta})=\check {\xi}$, then
$p\restriction P_{\mu}\Vdash \underline{f}(\check{\zeta})=\check {\xi}$. 
\end{claim}

\begin{proof}[Proof of the Claim]
    Assume on the contrary that 
    $r_0\le p\restriction P_{\mu}$ and $r_0\Vdash \underline{f}(\check {\zeta})\ne \check {\xi}$.
    Let $r_1=(r_0\restriction P_{\mu})\cup (p\setm  (p\restriction P_{\mu}))$. Then $r_1\le p$ and so 
    $r_1 \Vdash\underline{c}(\check {\zeta})=\check {\xi}$.
Since $r_0\restriction P_{\mu}=r_1\restriction P_{\mu}$,
    there is ${\pi}\in H_{\mu}$ such that $r_0$ and ${\pi}(r_1)$
    are compatible. 
    Then $$r_0\land {\pi}(r_1)\Vdash \underline{f}(\check {\zeta})=\check {\xi}
    \land \underline{f}(\check {\zeta})\ne\check{\xi}.$$
    Contradiction. We proved the Claim.
\end{proof}
By the Claim, if we take 
\begin{displaymath}
\underline{h}=\{\big\langle \widecheck{\<{\zeta},{\xi}\>}^{P_{\mu}},p\big\rangle:p\in P_{\mu}\land p\Vdash \underline{f}(\check {\zeta})=\check {\xi} \}
\end{displaymath}
then $\KK{\underline{h}}{\mc G\restriction P_{\mu}} =f$. So
we proved the Key lemma.
\end{proof}

\begin{lemma}\label{lm:zfcsubmodels}
For each ${\mu}<{\kappa}$, $\mc M[\mc G\restriction P_{\mu}]\subs \mc N_1$.
\end{lemma}

\begin{proof}
Indeed, if $\underline{a}$ is a $P_{\mu}$-name, then $st(\underline{a})\ge H_{\mu}$.
\end{proof}

\begin{lemma}
    For each ${\lambda}\in C$, $\mc N\models$ ${\lambda}$ is countable.  
    \end{lemma}
    
\begin{proof}
We  have  $\mc M[\mc G\restriction P_{\mu+1}]\models$ "${\mu}$  is countable"
and   $\mc M[\mc G\restriction P_{\mu+1}]\subs \mc N_1$.
\end{proof}

\begin{lemma}
    $\mc N\models$ $cf({\omega}_1)={\omega}_1={\kappa}$.   
   \end{lemma}
   
   \begin{proof}
       $P$ satisfies ${\kappa}$-cc.  So 
       $\mc M[\mc G]\models $ "${\kappa}={\omega}_1$" and so ${\kappa}$ is regular in $\mc N_1\subs \mc M[\mc G]$. 
Thus, ${\kappa}$ is a regular cardinal in $\mc N_1$. So ${\kappa}={\omega}_1$ in $\mc N_1$.
    \end{proof}

\newpage

\begin{lemma}
    $\mc N_1\models $ ${\omega}_1\to ({\nu})^n_{\omega}$ for each 
    ${\nu}<{\omega}_1$.
    \end{lemma}


\begin{proof}
   By the key lemma, there is ${\mu}<{\kappa}$ with $c\in \mc M[\mc G\restriction P_{\mu}]$.

Since ${\kappa}$ is inaccessible in  $\mc M[\mc G\restriction P_{\mu}]$, we can  
fix   ${\lambda}\in C$ with
\begin{displaymath}
    \mc M[\mc G\restriction P_{\mu}]\models
    \exp_{n-1}({\nu})<{\lambda}. 
\end{displaymath}   

In $\mc M[\mc G\restriction P_{\mu}]$ consider the coloring 
$c\restriction {[{\lambda}]}^{n}:{[{\lambda}]}^{n}\to {\omega}$.

By the Erdős-Rado theorem $\exp_{n-1}({\nu})^+\to ({\nu}^+)^n_{\nu}$, 
\begin{multline*}
    \mc M[\mc G\restriction P_{\mu}]\models\text{there is set $A\subs {\lambda}$
    with order type ${\nu}$}\\\text{ 
    which is $c$-homogeneous. 
    } 
\end{multline*}
Then $A \in \mc N_1$, $A$ has order type ${\nu}$ and $A$ is $c$-homogeneous in 
$\mc N_1$ as well. 
\end{proof}


\begin{lemma}\label{lm:wc}
    If ${\kappa}$ is weakly compact, then 
    \begin{displaymath}
    \mc N_1\models {\omega}_1\to ({\omega}_1)^n_{\omega}
    \end{displaymath}
    for each $n\in {\omega}$.
    \end{lemma}
        
    \begin{proof}
    Assume that 
    \begin{displaymath}
    \mc N_1\models f:{[{\omega}_1]}^{n}\to {\omega}.
    \end{displaymath}
    By the Key Lemma \ref{lm:key} there is ${\mu}<{\kappa}$ with 
    $f\in \mc M[\mc G\cap P_{\mu}]$. The cardinal ${\kappa}$ is weakly compact in $\mc M[\mc G\cap P_{\mu}]$,
    so there is an $f$-homogeneous $A\in {[{\kappa}]}^{{\kappa}}\cap \mc M[\mc G\cap P_{\mu}]$.
    Hence, $A$ has a name $\underline A$ with $st(\underline A)\ge H_{\mu}$.
    Thus,  $A\in \mc N_1$ and 
    \begin{displaymath}
    \mc N_1\models \text{$A$ is $f$-homogeneous with cardinality ${\omega}_1$.}
    \end{displaymath}
    \end{proof}

\begin{lemma}\label{lm:m1ed}
    $\mc N_1\models $ $M({\omega}_1,{\omega},2)\to ED$ fails.
    \end{lemma}

    \begin{proof}
        For ${\alpha}<{\kappa}$ let $$I_{\alpha}=\{{\alpha}+n:n<{\omega}\},$$
for ${\lambda}\in C$ and $n<{\omega}$ let 
\begin{displaymath}
D_{{\lambda},n}=\{{\lambda}+{\omega}{\alpha}+n: {\alpha}<{\lambda}\}.
\end{displaymath}
and let 
\begin{displaymath}
\mc A=\{I_{\alpha}:{\alpha}<{\kappa}\text{ is limit}\}\cup\{D_{{\lambda},n}:{\lambda}\in C,n<{\omega}\}.
\end{displaymath}
Then $\mc A\in \mc M$. Since $\mc A\subs {[{\kappa}]}^{<{\kappa}}$ in $\mc M$,
we have  $$\mc N_1\models  \mc A\subs {[{\omega}_1]}^{{\omega}}.$$ 
Clearly $\mc A$ is 2-almost disjoint.

For each $A\in \mc A$ let $F(A)$ be the first $2$ elements of $A$ in the natural orderings 
of ${\kappa}$. Since the assignment  $A\mapsto F(A)$ is injective, and 
${[{\kappa}]}^{2}$ has a well-ordering in type ${\kappa}$, 
we can assume that there is bijection  $A$ 
between ${\kappa}$ and $\mc A$.

Define $B:{\kappa}\times {\omega}_1\to {\kappa}$ as follows:
$B({\alpha},{\zeta})$ is the ${\xi}$th element of $A({\alpha})$ provided the 
order type of $A({\alpha})$ is greater than ${\xi}$, and 
$B({\alpha},{\zeta})=\min A({\alpha})$ otherwise. 

Assume that $F$ witnesses that  $\mc A$ is ED.
We can assume that $F(A)\ne \empt$ for each $A\in \mc A$.
Define $E:{\kappa}\times {\omega}\to {\kappa}$ as follows: 
$E({\alpha},n)$ is the $n$th element of $F(A({\alpha}))$ provided that
 $|F(A({\alpha}))|\ge n$, and 
$E({\alpha},n)=\min F(A({\alpha}))$ otherwise.

By the key lemma, there is ${\mu}<{\kappa}$ such that 
$B,E\in \mc M[\mc G\restriction P_{\mu}]$.
Then $F,\mc A\in \mc M[\mc G\restriction P_{\mu}]$

Thus, $\mc A$ is ED in $\mc M[\mc G\restriction P_{\mu}]$.
But in that model $({\omega}_1)^{\mc M[\mc G\restriction P_{\mu}]}={\lambda}$
for some ${\lambda}<{\kappa}$.

Let $H=\bigcup\{F(D_{{\lambda},n}:n<{\omega}\}$.
Since $H$ is a countable union of finite sets, $|H|\le {\omega}$ by Theorem 
\ref{tm:elementary-in-ZF}(\ref{en:cof-zf}).
Thus, there is ${\alpha}<{\lambda}$ such that $I_{{\lambda}+{\alpha}}\cap H=\empt$.

Let $n\in {\omega}$ such that ${\lambda}+{\omega}{\alpha}+n\notin F(I_{{\lambda}+{\alpha}})$.
Thus, 
$${\lambda}+{\omega}{\alpha}+n\in 
\big (D_{{\lambda},n}\setm F(D_{{\lambda},n}\big)
\cap \big(I_{{\lambda}+{\alpha}}\setm F(I_{{\lambda}+{\alpha}})\big).$$
Contradiction.
\end{proof}

\begin{lemma}
    $\mc N_1\models$   There is no uniform denumeration on ${\omega}_1$. 
\end{lemma}

\begin{proof}
Put together Lemma \ref{lm:m1ed} and Theorem \ref{tm:ud2ed-gen}.
%
\end{proof}

So we proved that (1) implies (2).

\noindent (2)$\to$(1)

We can proof a bit stronger statement:

\begin{claim}\label{cl :stronger}
    If ZF + $cf({\omega}_1)={\omega}_1$ + $\neg$UD$({\omega}_1)$  is consistent then so is ZFC + there is an inaccessible cardinal. 
\end{claim}


Since ${\omega}_1$ is regular, $L\models \text{ ``${{\omega}_1}^V$ is a regular cardinal''}$. 
If $$L\models \text{ ``${{\omega}_1}^V={\kappa}^+$ for some cardinal ${\kappa}\in On $'' }$$
then there is a function $f\in L$ with $\dom(f)={{\omega}_1}^V= ({\kappa}^+)^L$ 
such that $f({\alpha})$ is a functions mapping  ${\kappa}$ onto ${\alpha}$ for 
$0<{\alpha}<({\kappa}^+){}^L$.
Define  
$$F=\{\<{\alpha},{\zeta}, f({\alpha})({\zeta})\>:0<{\alpha}<({\kappa}^+)^L,{\zeta}<{\kappa}\},$$
and let $g\in V$ be a bijection between ${\omega}$ and ${\kappa}$.
Then $f,g\in  L[F,g]$ and 
\begin{displaymath}
L[F,g]\models (f({\beta})\circ g)[{\omega}]={\beta}\text{ for each }  {\beta}\in {{\omega}_1}^V\setm {\kappa}.
\end{displaymath}
Thus, ${\omega}_1^{L[F,g]}={\omega}_1$, and so UD$({\omega}_1)$ holds by lemma \ref{lm:un-eq},
which is not the case. 
So ${{\omega}_1}^V$ is an inaccessible cardinal in $L$.
\end{proof}

\begin{problem}[ZF]\label{pr:ud2ed}
    (1)  Does $M({\omega}_1,  {\omega},n)\to ED$   imply  UD(${\omega}_1$)?
    (2)    Does ZF + UD(${\omega}_1$) imply  
    $M({\kappa},  {\omega},n)\to ED$ for each infinite cardinal ${\kappa}\in On $ 
    and natural number  $n<{\omega}$?

\end{problem}

\bibliographystyle{abbrv} 
       
\bibliography{inf-comb-zf}

\end{document}